\input amstex
\documentstyle{amsppt}
\NoBlackBoxes
\define\BK{1}
\define\G{2}
\define\Ha{3}
\define\Hb{4}
\define\LM{5}
\define\Li{6}
\define\Mc{7}
\define\Z{8}
\define\ga{\gamma}
\define\la{\lambda}
\define\zt{\zeta}
\define\Sym{\operatorname{Sym}}
\define\HH{\frak H}
\topmatter
\title Periods of mirrors and multiple zeta values\endtitle
\author Michael E. Hoffman \endauthor
\address U. S. Naval Academy, Annapolis, MD 21402\endaddress
\email meh\@nadn.navy.mil\endemail
\abstract
In a recent paper, A. Libgober showed that the multiplicative sequence
$\{Q_i(c_1,\dots,c_i)\}$ of Chern classes corresponding to the 
power series $Q(z)=\Gamma(1+z)^{-1}$ appears in a relation between
the Chern classes of certain Calabi-Yau manifolds and the periods of their
mirrors.  We show that the polynomials $Q_i$ can be expressed in terms
of multiple zeta values.
\endabstract
\endtopmatter
\document
\head 1. The multiplicative sequence\endhead
In \cite{\Li}, the (Hirzebruch) multiplicative sequence $\{Q_i\}$ associated to
the power series $Q(z)=\Gamma(1+z)^{-1}$ is considered in connection with
mirror symmetry.
If $e_i$ denotes the $i$th elementary symmetric function in the variables 
$t_1, t_2,\dots$, then
$$
\sum_{i=0}^{\infty}Q_i(e_1,\dots,e_i)=\prod_{i=1}^{\infty}\frac1{\Gamma 
(1+t_i)} .
\tag1$$
As shown in \cite{\Li}, the polynomials $Q_i(c_1,\dots,c_i)$ in the
Chern classes of certain Calabi-Yau manifolds $X$ are related to
the coefficients of the generalized hypergeometric series expansion of 
the period (holomorphic at a maximum degeneracy point) of a mirror of $X$.  
In particular, if $X$ is a Calabi-Yau hypersurface
of dimension 4 in a nonsingular toric Fano manifold, then
$$
\int_X Q_4(c_1,c_2,c_3,c_4)=\frac1{24}K_{ijkl}\frac{\partial^4 c(0,\dots,0)}
{\partial\rho_i \partial\rho_j \partial\rho_k \partial\rho_l} ,
$$
where the $c(\rho_1,\dots,\rho_r)$ are coefficients in the expansion of
the period and $K_{ijkl}$ is the (suitably normalized)
4-point function corresponding to a mirror of $X$.  In \cite{\Li} it is 
shown that the
polynomials $Q_i$ have the form 
$$
Q_1(c_1) = \gamma c_1\quad\text{and}\quad
Q_i(c_1,\dots,c_i)=\zt(i) c_i + \cdots , i>1.
$$
In this note we show that the polynomials $Q_i$ have an explicit 
expression in terms of multiple zeta values (called multiple harmonic
series in \cite{\Hb}), which have previously appeared in connection
with Kontsevich's invariant in knot theory \cite{\Z,\LM}, and in 
quantum field theory \cite{\BK}.
\head 2. The formula for the $Q_i$\endhead
Let Sym be the algebra of symmetric functions in the variables
$t_1,t_2,\dots$ (with rational coefficients), and let $p_i$ be the 
$i$th power-sum symmetric function in these variables.
For a partition $\la=(\la_1,\la_2,\dots)$, let $m_{\la}$ be the
corresponding monomial symmetric function and $e_{\la}=e_{\la_1}
e_{\la_2}\cdots$.  It is well known that $\{m_{\la}\}$ and $\{e_{\la}\}$
are bases for $\Sym$ as a vector space.
In \cite{\Hb} it is shown (Theorem 5.1) that the
homomorphism $\zt:\Sym\to\bold R$ such that $\zt(p_1)=\gamma$
and $\zt(p_i)=\zt(i)$ for $i>2$ satisfies
$$
\zt\left(\sum_{i\ge 0} e_i z^i\right)=\frac1{\Gamma(1+z)} .
\tag2$$
Our main result expresses the polynomials $Q_i$ in terms of $\zt$.
\proclaim{Theorem} For any partition $\la$ of $i$, the coefficient
of $e_{\la}$ in $Q_i(e_1,\dots,e_i)$ is $\zt(m_{\la})$.
\endproclaim
\demo{Proof} Using equations (1) and (2), we have
$$
\sum_{i\ge 0}Q_i(e_1,e_2,\dots)=
\prod_{i=1}^{\infty}\frac1{\Gamma(1+t_i)}=
\prod_{i=1}^{\infty}\sum_{j=0}^{\infty}\zt(e_j)t_i^j =
\sum_{\la}\zt(e_{\la})m_{\la} .
$$
Now the transition matrix $M$ from the basis $\{e_{\la}\}$ of 
$\Sym$ to the basis $\{m_{\la}\}$, i.e.
$$
e_{\la}=\sum_{\mu}M_{\la \mu}m_{\mu}
$$
is known to be symmetric (see Ch. I, \S 6 of \cite{\Mc}), so we have
$$
\sum_{\la}\zt(e_{\la})m_{\la}=
\sum_{\la}\sum_{\mu}M_{\la \mu}\zt(m_{\mu})m_{\la}=
\sum_{\mu}\zt(m_{\mu})\sum_{\la}M_{\mu \la}m_{\la}=
\sum_{\mu}\zt(m_{\mu})e_{\mu} ,
$$
and the result follows.
\enddemo
\head 3. Multiple zeta values\endhead As shown in \cite{\Hb}, $\zt$
can be thought of as a homomorphism from the algebra of quasi-symmetric 
functions in the $t_i$ (as defined in \cite{\G}) to $\bold R$ that extends 
the multiple zeta values introduced in \cite{\Ha} and \cite{\Z}, i.e.
$$
\zt(i_1,i_2,\dots,i_k)=\sum_{n_1>n_2>\dots>n_k\ge1}\frac1{n_1^{i_1}
n_2^{i_2}\cdots n_k^{i_k}} ,
\tag3$$
where $i_1$ must be assumed greater than 1 for convergence.
If we let $\HH^1$ be the rational vector space of polynomials in
the noncommuting
variables $z_1,z_2,\dots$, then $\HH^1$ becomes isomorphic to the
algebra of quasi-symmetric functions if we define the 
(commutative) multiplication * by the inductive rule
$$
z_i w_1 * z_j w_2 = z_i(w_1*z_j w_2)+z_j(z_i w_1 * w_2)+z_{i+j}(w_1*w_2)
$$
for any words $w_1, w_2$ in the $z_i$; see \cite{\Hb} for details.
The algebra Sym of symmetric functions can be identified with the subspace
of $\HH^1$ generated by linear combinations of monomials invariant under
permutation of subscripts, e.g. $z_2^2=m_{22}$ and $z_1z_2+z_2z_1=m_{21}$;
$z_i$ and $z_1^i$ correspond to $p_i$ and $e_i$ respectively.
As an algebra $\HH^1$ is generated by Lyndon words in the $z_i$, i.e.
monomials $w$ such that for any nontrivial decomposition $w = uv$ one
has $v>w$, where the $z_i$ are ordered as $z_1>z_2>\cdots$ and this order
is extended to monomials lexicographically.  Then the only Lyndon word
that starts with $z_1$ is $z_1$ itself, and $\zt$ is the homomorphism
from $\HH^1$ to $\bold R$ defined on Lyndon words $w=z_{i_1}z_{i_2}\cdots
z_{i_k}$ by
$$
\zt(w)=\cases \gamma,&w=z_1,\\ \zt(i_1,i_2,\dots,i_k),&\text{otherwise.}
\endcases
$$
By the results of \cite{\Hb}, $\zt(z_{i_1}z_{i_2}\cdots z_{i_k})$ coincides
with $\zt(i_1,i_2,\dots,i_k)$ as defined by equation (3) whenever $i_1>1$.
\par
Since the power-sum symmetric functions $p_i$ generate the algebra
$\Sym$, we can compute $\zt(m_{\la})$ by first expressing $m_{\la}$ in
terms of power-sum functions (see \cite{\Mc, p. 109} for an explicit
formula), and then applying the homomorphism $\zt$.  Hence the coefficient
of each monomial $c_{\la}$ in $Q_i(c_1,\dots,c_i)$ is a polynomial in
the numbers $\ga$ and $\zt(i)$, $i\ge 2$.  By this method we can obtain
formulas (1.3)-(1.6) of \cite{\Li} (with the following corrections:
in (1.4) the coefficient of $c_1^2$ should be $\frac12(\ga^2-\zt(2))$,
while in (1.5) the coefficient of $c_1^3$ should be $\frac13\zt(3)-\frac
12\ga\zt(2)+\frac16\ga^3$).
\par
If $m_{\la}$ is a monomial symmetric function such that the partition
$\la$ involves no 1's, then $\zt(m_{\la})$ is just a sum of ordinary multiple 
zeta values, e.g. $\zt(m_{22})=\zt(2,2)=\frac34\zt(4)$ and $\zt(m_{62})=
\zt(6,2)+\zt(2,6)=\frac23\zt(8)$.
In particular, such $\zt(m_{\la})$ are the only coefficients
needed to evaluate $Q_i(c_1,\dots,c_i)$ 
on a Calabi-Yau manifold, since $c_1=0$ in that case.

\enddocument